\documentclass[11pt, draft]{amsart}
\usepackage{amssymb}
\usepackage[mathscr]{eucal}

\theoremstyle{plain}
\newtheorem{theorem}{Theorem}

\theoremstyle{remark}
\newtheorem*{remark}{Remark}

\newcommand\la{\langle}
\newcommand\ot{\otimes}
\newcommand\ra{\rangle}
\newcommand\rank{\operatorname{rank}}
\newcommand\tr{\operatorname{tr}}
\newcommand\T{\mathbb T}

\begin{document}
\title[]{On certain automorphisms of sets of partial isometries}
\author{LAJOS MOLN\'AR}
\address{Institute of Mathematics and Informatics\\
         University of Debrecen\\
         4010 Debrecen, P.O.Box 12, Hungary}
\email{molnarl@math.klte.hu}
\thanks{  This research was supported by the
          Hungarian National Foundation for Scientific Research
          (OTKA), Grant No. T030082, T031995, and by
          the Ministry of Education, Hungary, Reg.
          No. FKFP 0349/2000}
\date{\today}
\begin{abstract}
Under the mild condition of continuity at a single point we describe all
the bijections
of the set of all partial isometries on a Hilbert space which preserve
the order and the orthogonality in both directions. Moreover,
we present a natural analogue of Wigner's theorem on quantum
mechanical symmetries for the set of all rank-1 partial isometries.
\end{abstract}
\maketitle

In the mathematical foundations of quantum mechanics the lattice
$P(H)$ of all projections on a Hilbert space $H$ plays a fundamental
role. The whole set of
projections equipped with the usual partial ordering and orthogonality
represents the probabilistic aspect of the theory and the set $P_1(H)$
of all
rank-1 projections with the notion of transition probability is the
object of the fundamental theorem of Wigner on quantum
mechanical symmetries.
Clearly, the projections can be characterized as positive partial
isometries. If we drop positivity, one can naturally raise
several problems
concerning partial isometries which are familiar in relation to
projections.
Probably the most fundamental results concerning the mentioned
structures of projections
are the description  of all bijective transformations of $P(H)$ which
preserve the order and the orthogonality in both directions, and
Wigner's
theorem (sometimes called unitary-antiunitary theorem) determining all
bijective transformations of $P_1(H)$ which preserve the transition
probabilities
(these results can be found, for example, in \cite{CVLL}; also see the
references therein).
The aim of this paper is to present analogue results for the
corresponding structures of partial isometries.
These operators, we mean the partial isometries, play the same role in
the relatively new theory of ternary structures of operators as the
projections do in connection with the usual binary structures. The
former structures
are of great importance in several subareas of functional analysis, for
example, in the study of the geometrical properties of operator algebras
(see \cite{DFR} and its references).

We begin with the notation that we use throughout.
Let $H$ be a complex separable Hilbert space.
The algebra of all bounded linear operators on $H$ is denoted by $B(H)$.
The symbol $F(H)$ stands for the set of all finite rank elements of
$B(H)$.
An operator $U\in B(H)$ is called a partial
isometry if there is a closed subspace $I_U$ of $H$ called the initial
space of $U$ for which it holds that $U$ acts as an isometry on $I_U$
and $U$ is zero on $I_U^\perp$. The range of $U$ is called the final
space of $U$ and it is denoted by $F_U$.
The set of all partial isometries on $H$ is denoted by $PI(H)$.
The symbol $PI_1(H)$ stands for the set of all rank-1 partial isometries
on $H$. We introduce a partial ordering on $PI(H)$. For any $P,Q\in
PI(H)$ we write $P\leq Q$ if $I_P\subset I_Q$, $F_P\subset F_Q$ and
$P_{|I_P}=Q_{|I_P}$. We say that the partial isometry $P$ is orthogonal
to the partial isometry $Q$ if $P^*Q=PQ^*=0$.
Obviously, the partial ordering introduced above
coincides, when restricting to projections, with the usual
order on $P(H)$.
The same holds true for the orthogonality
relation.

A linear norm preserving bijection on $H$ is called a unitary operator ,
while any conjugate-linear norm preserving bijection on $H$ is called an
antiunitary operator.

If $x, y\in H$, then $x\ot y$ stands for the operator defined by
\[
(x\ot y)(z) = \la z,y\ra x \qquad (z\in H).
\]

The usual trace functional on the trace class operators is denoted by
$\tr$.
Finally, the set of all complex numbers of modulus 1 is denoted by
$\mathbb T$.

Our first result which follows describes all the bijections of $PI(H)$
which preserve the order and orthogonality in both directions.

\begin{theorem}
Let $H$ be a complex Hilbert space with $\dim H \geq 3$. Suppose that
$\phi :PI(H) \to PI(H)$ is a bijective transformation which preserves
the
partial ordering and the orthogonality between partial isometries
in both directions.
If $\phi$ is continuous (in the operator norm)
at a single element of $PI(H)$ different from 0, then $\phi$ can be
written in one of the following forms:
\begin{itemize}
\item[(i)]
there exist unitaries $U,V$ on $H$ such that
\[
\phi(R)=URV \qquad (R\in PI(H));
\]
\item[(ii)]
there exist antiunitaries $U,V$ on $H$ such that
\[
\phi(R)=URV \qquad (R\in PI(H));
\]
\item[(iii)]
there exist unitaries $U,V$ on $H$ such that
\[
\phi(R)=UR^*V \qquad (R\in PI(H));
\]
\item[(iv)]
there exist antiunitaries $U,V$ on $H$ such that
\[
\phi(R)=UR^*V \qquad (R\in PI(H)).
\]
\end{itemize}
\end{theorem}

\begin{proof}
We assert that $P\in PI(H)$ is a rank-$n$ partial isometry if and only
if so is $\phi(P)$. By the order preserving property of $\phi$ we
readily have $\phi(0)=0$ which verifies our claim in the case
$n=0$.
Suppose that we have the assertion for $k=0, \ldots, n$. Let $P$ be a
partial isometry on $H$. It is clear that $P$ has rank $n+1$ if
and only if for every $Q\in PI(H)$ with $Q\neq P$, $Q\leq P$ it follows
that the rank of $Q$ is less than or equal to $n$ and there is such a
$Q$ whose rank is $n$. The order preserving
property of $\phi$ now implies that the rank of $\phi(P)$ is $n+1$.

Our next aim is to show that $\phi$ is completely orthoadditive.
The meaning of this property will be clear in a moment.
Let $R_\alpha \in PI(H)$ be a set of pairwise orthogonal partial
isometries. It is well-known that the series
$\sum_\alpha R_\alpha$ is convergent in the strong operator topology and
its sum $R$ is a partial isometry. The sum $\sum_\alpha
\phi(R_\alpha)$ is also a partial isometry, so there exists a $Q\in
PI(H)$ such that
\[
\sum_\alpha \phi(R_\alpha)=\phi(Q).
\]
Since $R_\alpha \leq R$ for every $\alpha$, it follows that
$\phi(Q)=\sum_\alpha \phi(R_\alpha)\leq \phi(R)$. On the other
hand, since $\phi(R_\alpha) \leq \phi(Q)$, we have $R_\alpha \leq Q$ for
every $\alpha$ and we obtain
that $R=\sum_\alpha R_\alpha \leq Q$. This yields $\phi(Q)=\phi(R)$.
Therefore, we have
\[
\sum_\alpha \phi(R_\alpha)=\phi(\sum_\alpha R_\alpha)
\]
and this means the complete orthoadditivity of $\phi$.

Take $3\leq n\in \mathbb N$.
Let $M, N$ be $n$-dimensional subspaces of $H$. Denote by $PI(M,N)$
the set of all partial isometries on $H$ whose initial space is a subset
of $M$ and whose final space is a subset of $N$. Let $P$ be a partial
isometry
on $H$ whose initial space is $M$ and whose final space is $N$. Then
$\phi(P)$ is a rank-$n$ partial isometry. Let $M'=I_{\phi(P)}$ and
$N'=F_{\phi(P)}$. It is easy to see that an operator $Q\in PI(H)$
belongs to $PI(M,N)$ if and only if for every $R\in PI(H)$ which is
orthogonal to $P$ it follows that $R$ is orthogonal to $Q$.
By the orthogonality preserving property of $\phi$ we deduce that
$\phi$ maps $PI(M,N)$ onto $PI(M',N')$.
Clearly, $PI(M,N)$ and $PI(M',N')$ are both isomorphic to
the space $PI(H_n)$ of all partial isometries on an $n$-dimensional
Hilbert space $H_n$.

Therefore, our map
$\phi$ induces a transformation $\psi$ on $PI(H_n)$
which has the same preserver properties as $\phi$.
Obviously, $\psi(I)$ is unitary. Multiplying $\psi$ by a
fixed unitary element of $B(H_n)$, we can assume that $\psi(I)=I$. This
latter property implies that $\psi$ sends projections to projections.
In fact, for any partial isometry $P$ we have $P\leq I$ if and only if
$P$ is a projection. So, $\psi$ is a bijective transformation
of the set of all projections in $B(H_n)$ which preserves the order and
the orthogonality in both directions.
The form of such transformations is well-known.
It follows, for example, from \cite{CVLL}
that there exists an either
unitary or antiunitary operator $U$ on $H_n$ such that
\begin{equation}\label{E:partial2}
\psi(P)=UPU^*
\end{equation}
for every projection $P$ in $B(H_n)$.

Let $\lambda\in \T$. We state that
$\psi(\lambda I)$ is a scalar multiple of the identity. To see this,
first observe that if
$R\in PI(H)$ is a rank-1 partial isometry, then $\phi(\lambda R)$ is a
scalar multiple of $\phi(R)$ (which scalar obviously has modulus 1).
This follows from the fact that
any partial isometry in $B(H)$ is orthogonal to $\phi(\lambda R)$ if and
only if it is orthogonal to $\phi(R)$ and that $\phi(\lambda R)$,
$\phi(R)$ are both of rank 1.
Now, one can characterize the scalar operators in $PI(H_n)$ in the
following way.
The partial isometry $A\in PI(H_n)$ is equal to the identity multiplied
by a scalar from $\T$ if and only if for every rank-one projection $P$
we have that $\lambda P \leq A$ for some $\lambda\in \T$.
Therefore, we obtain that $\psi$
preserves the scalar partial isometries.
Consequently, there is a function $f: \T \to \T$ such that
\[
\psi(\lambda I)=f(\lambda) I.
\]

We know that
if $P\in B(H_n)$ is any rank-one projection and $\lambda \in \T$, then
there is a $\mu\in \T$ such that $\mu \psi(P)= \psi(\lambda P)$. As
$\psi(\lambda P)\leq \psi(\lambda I)$, we obtain that $\mu \psi(P)\leq
f(\lambda) I$. This implies that $\mu=f(\lambda)$, so we have
\[
\psi(\lambda P)=f(\lambda )\psi(P)
\]
for every rank-1 projection $P\in B(H_n)$ and $\lambda \in \T$.
By the orthoadditivity of $\psi$ (this follows form the
orthoadditivity of $\phi$) we infer that
the previous equality holds true for every projection $P\in B(H_n)$
without any restriction on its rank.

If $R\in PI(H_n)$ is any partial isometry, then $R$ can be written in
the form $R=UP$, where $U$ is unitary and $P$ is a projection.
Considering the
transformation $Q \mapsto \psi(U)^*\psi(UQ)$ on $PI(H_n)$
and applying what we have proved above, it follows that there is a
function $f_R: \T \to \T$ such that
\[
\psi(\lambda R)=f_R(\lambda )\psi(R) \qquad (\lambda \in \T).
\]
The function $f_R$ might depend on $R$ since the function $f$ appearing
above depends on $\psi$. However, we prove that $f_R=f$ for
every partial isometry $R\in PI(H_n)$. This will be done below.

Let $R\in PI(H_n)$ be a rank-1 partial isometry. Pick a rank-1
projection $P\in B(H_n)$ which is orthogonal to $R$.
By the orthoadditivity of $\psi$, for any $\lambda \in \T$ we compute
\[
\psi(\lambda R+\lambda P)=
f_{R+P}(\lambda )\psi(R+P)=
f_{R+P}(\lambda )(\psi(R)+\psi(P)).
\]
On the other hand, we have
\[
\psi(\lambda R+\lambda P)=
\psi(\lambda R)+\psi(\lambda P)=
f_R(\lambda)\psi(R)+f(\lambda)\psi(P).
\]
These imply that $f_{R}=f_{R+P}=f$. So we have
\begin{equation}\label{E:partial1}
\psi(\lambda R)=f(\lambda) \psi(R) \qquad (\lambda \in \T)
\end{equation}
for every rank-1 partial isometry $R\in PI(H_n)$. Since every partial
isometry is the sum of mutually orthogonal rank-1 partial isometries,
by the orthoadditivity of $\psi$ we obtain
that the above equality holds also for every partial isometry $R$ on
$H_n$.

Since for every $\lambda, \mu \in \T$ we have
\[
f(\lambda \mu) I=\psi((\lambda \mu) I)=\psi(\lambda (\mu I))=f(\lambda)
f(\mu)I,
\]
it follows that $f:\T \to \T$ is a multiplicative bijection.

We supposed that our original transformation is norm-continuous at a
point $0\neq R\in PI(H)$. We can find pairwise
orthogonal rank-1 partial isometries $R_\alpha$ for which $\sum_\alpha
R_\alpha=R$. By the complete orthoadditivity of $\phi$ we obtain that
\[
\sum_\alpha \phi(\lambda R_\alpha) =\phi(\lambda R).
\]
As we have seen above, the operator $\phi(\lambda R_\alpha)$ is a scalar
multiple of $\phi(R_\alpha)$ for every $\lambda\in \T$.
Let $\lambda_n$ be a sequence in $\T$ converging to 1. This implies that
$\phi(\lambda_n R) \to \phi(R)$ and then we have
$\phi(\lambda_n R_\alpha) \overset{n\to \infty}{ \longrightarrow }
\phi(R_\alpha)$ for every $\alpha$.
Pick an $\alpha_0$ and let us suppose that $R_{\alpha_0} \in PI(M,N)$.
We state that $f$ is continuous.
In fact, by \eqref{E:partial1} we can deduce that
$f$ is continuous at 1. Since $f$ is multiplicative, we
obtain that $f$ is a continuous character of $\T$. It is well-known that
the continuous characters of $\T$ are exactly the functions $z
\mapsto z^n$ $(n \in \mathbb Z)$.
Since $f$ is bijective, we obtain that $f$ is either the
identity or the conjugation on $\T$.

Suppose that the operator $U$ in \eqref{E:partial2} is unitary and
suppose that $f$ is the conjugation on $\T$.
Using the spectral theorem and the orthoadditivity of $\psi$, we obtain
from \eqref{E:partial1} that for every unitary $R$ in $B(H_n)$ we have
\[
\psi(R)=UR^*U^*.
\]
If $R\in PI(H_n)$ is arbitrary, then there is a partial
isometry $R'$ on $H_n$ which is orthogonal to $R$ such that $R+\lambda
R'$ is unitary for every $\lambda \in \T$. It follows that
\[
\psi(R)\leq
\psi(R+\lambda R')=
U(R+\lambda R')^*U^*=
UR^*U^*+{\overline{\lambda}} U{R'}^*U^*
\]
for every $\lambda\in \T$. This obviously implies that $\psi(R)\leq
UR^*U^*$. Since $\rank \psi(R)=\rank R= \rank UR^*U^*$,
we obtain $\psi(R)= UR^*U^*$.
Examining the remaining cases concerning $U$ and $f$, we find that
$\psi$ is of one of the following forms:
\begin{itemize}
\item[(i1)]
there exists a unitary $U$ on $H_n$ such that
\[
\psi(R)=URU^* \qquad (R\in PI(H_n));
\]
\item[(ii1)]
there exists a unitary $U$ on $H_n$ such that
\[
\psi(R)=UR^*U^* \qquad (R\in PI(H_n));
\]
\item[(iii1)]
there exists an antiunitary $U$ on $H_n$ such that
\[
\psi(R)=URU^* \qquad (R\in PI(H_n));
\]
\item[(iv1)]
there exists an antiunitary $U$ on $H_n$ such that
\[
\psi(R)=UR^*U^* \qquad (R\in PI(H_n)).
\]
\end{itemize}

Going back to our original map $\phi$, we see that on $PI(M,N)$ the
function $\phi$ is of one of the following forms:
\begin{itemize}
\item[(i2)]
there exist unitaries $U,V$ on $H$ such that
\[
\phi(R)=URV \qquad (R\in PI(M,N));
\]
\item[(ii2)]
there exist unitaries $U,V$ on $H$ such that
\[
\phi(R)=UR^*V \qquad (R\in PI(M,N));
\]
\item[(iii2)]
there exist antiunitaries $U,V$ on $H$ such that
\[
\phi(R)=URV \qquad (R\in PI(M,N));
\]
\item[(iv2)]
there exist antiunitaries $U,V$ on $H$ such that
\[
\phi(R)=UR^*V \qquad (R\in PI(M,N)).
\]
\end{itemize}

Pick two mutually orthogonal unit vectors $x,y\in H$ and consider the
partial isometries $R_1=x\ot x$, $R_2=i x\ot x$, $R_3=x\ot y$.
Denote by $[x,y]$ the subspace generated by $x,y$. On $PI([x,y], [x,y])$
$\phi$ is of one of the above forms. Suppose that we have
(ii2). It follows that $\phi(R_2)=-i\phi(R_1)$, $\phi(R_3)^*\phi(R_1)=
0$ and $\phi(R_3)\phi(R_1)^*\neq 0$.
Consider a finite
dimensional subspace $H_0$ of $H$ which contains all the initial and
final spaces of
$R_{\alpha_0},R_1,R_2,R_3$. It is easy to see that on
$PI(H_0, H_0)$ the transformation $\phi$ must be of the form (ii2). In
fact,
the other possibilities can be easily excluded considering the relations
among $\phi(R_1), \phi(R_2)$ and $\phi(R_3)$.

Let us take any finite rank operator $A\in F(H)$. By spectral theorem
and polar decomposition, $A$ can be written as a finite sum
\begin{equation}\label{E:partial4}
A=\sum_k \lambda_k P_k,
\end{equation}
where $P_k$'s are finite rank partial isometries and $\lambda_k$'s are
scalar. We define
\[
\Phi(A)=\sum_k {\overline{\lambda_k}} \phi(P_k).
\]
We claim that $\Phi$ is a conjugate-linear transformation on $F(H)$
extending
the restriction of $\phi$ onto the set of all finite rank partial
isometries on $H$. First, we have to show that $\Phi$ is well-defined.
Let
\[
A=\sum_l \mu_l Q_l
\]
be another resolution of $A$ similar to what was given in
\eqref{E:partial4}. Let $H_0\subset H_1$ be a finite dimensional
subspace
of $H$ which contains the initial and finial spaces of all $P_k,
Q_l$ appearing above. It follows that on
$PI(H_1, H_1)$  the transformation $\phi$ is of the form (ii2).
Therefore, we can write
\[
\sum_k {\overline{\lambda_k}} \phi(P_k) =UA^*V =
\sum_l {\overline{\mu_l}} \phi(Q_l)
\]
which means that $\Phi$ is well-defined. Now, the additivity and
conjugate-linearity of $\Phi$ is trivial to verify. So, we have a
conjugate-linear
map $\Phi: F(H) \to F(H)$ which extends the restriction of $\phi$ onto
the set of all finite rank partial isometries.
Using the 'local' form (ii2) of $\phi$, it is
obvious that
$\Phi$ satisfies $\Phi(AB^*A)=\Phi(A)\Phi(B)^*\Phi(A)$ $(A,B\in F(H))$,
that is $\Phi$ is a conjugate-linear Jordan-triple homomorphism.
Since every finite rank partial isometry is in the range of $\Phi$, it
follows that $\Phi$ is surjective. If $A\neq 0$, then
we can write $A$ in the form \eqref{E:partial4} where $P_k$'s are
pairwise orthogonal finite rank partial isometries. By the orthogonality
preserving property of $\phi$
we readily have $\Phi(A)\neq 0$. Therefore,
the transformation $A\mapsto \Phi(A^*)$
is a linear Jordan-triple automorphism of $F(H)$. Now, we can apply the
result \cite[Theorem 3]{MolnarSemrl} on the form of linear Jordan-triple
isomorphisms between standard operator algebras (in that paper we used
the term 'triple isomorphism'). This gives us that, on the set of all
finite
rank partial isometries, $\phi$ is one of the forms (ii2),
(iii2). Using the complete
orthoadditivity of $\phi$ we find that the same holds true
on the whole set $PI(H)$.

The last part of the proof has begun with supposing that $\phi$ is of
the form (ii2) on a certain subset $PI(H_0,H_0)$ of the set of all
partial isometries. In any of the other cases one can apply
very similar arguments.
\end{proof}

\begin{remark}
Observe that 0 is an isolated point of $PI(H)$, so the continuity of
$\phi$ at 0 would mean nothing.
We remark that it would be a nice improvement to show (if it is true at
all) that the result above holds without our assumption on
continuity.

Finally, we recall that it is quite common to refer to the analogue
result
concerning the lattice of projections as the fundamental theorem of
projective geometry.
\end{remark}

Wigner's celebrated unitary-antiunitary theorem can be formulated in
several ways. One possibility is the following (see, for example,
\cite{CVLL}). If $\phi:
P_1(H) \to P_1(H)$ is a bijective function which preserves the
transition probabilities, that is,
\[
\tr \phi(P)\phi(Q) =\tr PQ \qquad (P,Q \in P_1(H)),
\]
then there is an either unitary or antiunitary operator $U$ on $H$ such
that $\phi$ is of the form
\[
\phi(P)=UPU^* \qquad (P,Q \in P_1(H)).
\]
Our next result gives an analogue assertion concerning the set of all
rank-1 partial isometries.
For other Wigner-type results on other structures we refer to our recent
papers \cite{Molnar1}, \cite{Molnar2}.

\begin{theorem}
Let $\phi:PI_1(H) \to PI_1(H)$ be a bijective function with the property
that
\begin{equation}\label{E:partial3}
\tr \phi(P)^*\phi(Q) =\tr P^*Q \qquad (P,Q \in PI_1(H)).
\end{equation}
Then $\phi$ is of one of the following forms:
\begin{itemize}
\item[(a)]
there exist unitaries $U,V$ on $H$ such that
\[
\phi(R)=URV \qquad (R\in PI_1(H));
\]
\item[(b)]
there exist antiunitaries $U, V$ on $H$ such that
\[
\phi(R)=UR^*V \qquad (R\in PI_1(H));
\]
\end{itemize}
\end{theorem}

\begin{proof}
As we have noted in the proof of our first theorem,
if $A\in F(H)$, then $A$ can be written as a finite sum $A=\sum_k
\lambda_k P_k$ where $P_k$'s are rank-one partial isometries and
$\lambda_k$'s are scalars.
We define
\[
\Phi(A)=\sum_k \lambda_k \phi(P_k).
\]
We claim that $\Phi(A)$ is well-defined. Indeed, if
$A=\sum_l \mu_l Q_k$ is another resolution of $A$ of the above kind,
then we compute
\begin{equation}
\begin{gathered}
\tr \phi(R)^*(\sum_k \lambda_k \phi(P_k))=
\sum_k \lambda_k \tr \phi(R)^*\phi(P_k)=
\\
\sum_k \lambda_k \tr R^*P_k=
\tr R^*(\sum_k \lambda_k P_k)=
\tr R^*A
\end{gathered}
\end{equation}
and, similarly,
\[
\tr \phi(R)^*(\sum_l \mu_l \phi(Q_l))=
\tr R^*A.
\]
Therefore, it follows that
\[
\tr \phi(R)^*(\sum_k \lambda_k \phi(P_k))=
\tr \phi(R)^*(\sum_l \mu_l \phi(Q_l))
\]
for every $R\in PI_1(H)$. Since $\phi$ is surjective, we
obtain that
$\sum_k \lambda_k \phi(P_k)=
\sum_l \mu_l \phi(Q_l)$.
So, $\Phi$ is well-defined. It is now obvious that $\Phi$ is a linear
transformation on $F(H)$. Since the range of $\Phi$ contains every
rank-1 partial isometry, it follows that $\Phi$ is surjective.
An operator $A\in F(H)$ has rank 1 if and only if it is a nonzero scalar
multiple of a partial isometry. It follows that $\Phi$ preserves the
rank-1 operators in both directions. The form of all surjective linear
(even additive) transformations on $F(H)$ which preserve the rank-1
operators
is known. Namely, it follows from \cite[Theorem 3.3]{OmladicSemrl} that
either there are bijective linear operators $U,V$ on $H$ such that
\[
\Phi(x\ot y)=Ux\ot Vy \qquad (x,y \in H)
\]
or there are bijective conjugate-linear operators $U,V$ on $H$ such that
\[
\Phi(x\ot y)=Uy\ot Vx \qquad (x,y \in H).
\]
Depending on the actual case, we easily obtain from \eqref{E:partial3}
that $U,V$ are scalar multiples of unitaries or antiunitaries.
It is now trivial to complete the proof.
\end{proof}

\bibliographystyle{amsplain}

\end{document}